\documentclass{amsart}
\usepackage{enumerate,xypic}
\newtheorem{theo}{Theorem}
\DeclareMathOperator{\str}{str}
\DeclareMathOperator{\tr}{tr}
\DeclareMathOperator{\Str}{Str}
\DeclareMathOperator{\Hom}{Hom}
\begin{document}
\title{On the cyclic homology of supermatrices}
\author{Paul A. Blaga}
\address{``Babe\c{s}-Bolyai'' University, Faculty of Mathematics and Computer Sciences, 1, Kog\u{a}lniceanu Street,
400084 Cluj-Napoca, Romania} \email{pablaga@cs.ubbcluj.ro}
\begin{abstract}
The aim of this note is to show that the generalized supertrace, constructed in another paper of
the author, inducing an isomorphism between the Hochschild homology of a superalgebra and that of
the superalgebra of square supermatrices of a given type over $A$, induces, also, an isomorphism
between the cyclic homologies of the two superalgebras.
\end{abstract}

 \maketitle
\section{Superalgebras and supermatrices}
Superalgebras (i.e. $\mathbb{Z}_2$-graded algebras) are a very important particular case of graded
algebras that play an important role in modern mathematics and theoretical physics. They are,
particularly, central objects, as basic algebraic machinery in the theory of supermanifolds, which
are objects similar to the ordinary manifolds, regarded as ringed spaces, but the sheaf of smooth
functions on the manifolds is replaced by a suitable sheaf of superalgebras. Their different kind
of homologies can be defined in a way quite similar to the ungraded case, but we have to take into
account the grading, which results in extra signs added in suitable places.

We are not going to say much about superalgebras, they are described in details, for instance,
in~\cite{manin} or in \cite{bruzzo}. Nevertheless, a couple of things have to be mentioned. First
of all, in the case of superalgebras, the grading is called \emph{parity} and the homogeneous
elements are called \emph{even}, if their parity is 0 and \emph{odd} if their parity is 1.

Let us assume, once and for all, that $A=A_0\oplus A_1$ is a superalgebra with unit over a
commutative ring $k$. Usually, its is assumed that the ring also has a $\mathbb{Z}_2$-grading (i.e.
it is a superring). In most cases of interest, however, this is not so (actually, \emph{it is}, but
the grading is trivial, everything is in degree zero).

If $a\in A_i$ is a homogeneous element, its parity will be denoted by $|a|$. All the tensor
products in this paper are taken over the ground ring $k$ and should be thought of in the graded
sense (see, for instance, \cite{bruzzo} or \cite{maclane}). Also, the general rules for the
computation of the degrees in tensor products are used (of course, we have in mind, always, that
all the computations are made in $\mathbb{Z}_2$).

Defining matrices in the superalgebra is a little bit tricky. Let us denote, first, by $\Pi A$ the
algebra $A$, with the parity of the elements reversed. Then we define a \emph{free} $A$-module of
rank $p|q$ (see \cite{manin}) over $A$ to be a (graded) module which is isomorphic to the module
$A^{p|q}\equiv A^p\oplus \left(\Pi A\right)^q$. Clearly, $A^{p|q}$ has  a system a system of $p+q$
generators ($p$ are odd, and $q$ are even, whence the notation).

We want a matrix to represent, as in the ungraded case, a morphism between two free graded modules,
when homogeneous basis are fixed in the two modules. It turns out (\cite{manin},\cite{bruzzo}) that
the solution is to define a \emph{matrix in standard format} (or a \emph{supermatrix}, for short)
of type $(m|n)\times(p|q)$ with entries in the superalgebra $A$, to be just a $(m+n)\times (p+q)$
matrix $M$ with entries in $A$, together with a block decomposition of the form
\begin{equation}\label{supermatr}
M=
\begin{pmatrix}
M_{11}&M_{12}\\
M_{21}&M_{22}
\end{pmatrix},
\end{equation}
where  $M_{11}, M_{12}, M_{21}$ and $M_{22}$ are, respectively, matrices of type $m\times p$,
$m\times q$, $n\times p$ and $n\times q$. Let us denote by $\mathcal{M}((m|n),(p|q),A)$ the set of
all supermatrices of a given format. They have an obvious structure of $A$-module. If we define the
parity of a matrix in such a way that a matrix be even if $M_{11}$ and $M_{22}$ have even entries
and the other two -- odd entries and odd -- if the converse is true, then this $A$-module becomes a
\emph{graded} module, which is isomorphic, as expected, to the graded $A$-module
$\Hom_{A}\left(A^{m|n},A^{p|q}\right)$.

In what follows, we shall be interested in the case of \emph{square supermatrices} of type
$(p|q)\times (p|q)$. We shall use, for the module of this supermatrices, the simplified notation
$\mathcal{M}_{p,q}(A)$. With respect to the usual matrix multiplication, this $A$-module becomes,
as one can see immediately, a $k$-superalgebra, that will be called in the sequel \emph{the
superalgebra of supermatrices over $A$}, which is our main concern in this paper.
\section{The cyclic homology of superalgebras}
The cyclic homology for superalgebras is defined in a manner completely similar to the case of
ordinary, ungraded, algebras (see, for that, \cite{loday} or \cite{seibt}). We define, first, the
chain groups as
\begin{equation*}
C_n(A)=\underbrace{A\otimes \dots \otimes A}_{n+1 \;\text{copies}},
\end{equation*}
where all the tensor products are taken over the ground ring $k$, and, of course, they are
considered in the graded sense. Next,  we define the \emph{Hochschild differential} in the
following way (\cite{kassel}). Consider the \emph{face maps} $d_i:C_n(A)\to C_{n-1}(A)$,
$i=1,\dots, n$ ($n\geq 1$), defined as
\begin{equation*}
d_i(a_0\otimes \dots \otimes a_n)=
\begin{cases}
a_0\otimes \dots\otimes a_ia_{i+1}\otimes\dots \otimes a_n&\text{if}\quad 0\leq i <n\\
(-1)^{|a_n|(|a_0|+\dots|a_{n-1}|)}a_na_0\otimes a_1\otimes \dots\otimes a_{n-1}&\text{if}\quad i=n
\end{cases}
\end{equation*}
and the \emph{degeneracy maps} $s_i:C^n(A)\to C_{n+1}(A)$, $i=0,\dots, n$, defined by
\begin{equation*}
s_i(a_0\otimes \dots\otimes a_n)=a_0\otimes \dots\otimes a_i\otimes 1\otimes a_{i+1}\otimes\dots
\otimes a_n,\qquad, i=0,\dots, n.
\end{equation*}
 Then the
Hochschild differential is the map $b_n:C_n(A)\to C_{n-1}(A)$, given by
\begin{equation*}
b_n=\sum_{i=0}^n(-1)^i d_i.
\end{equation*}
We denote, also,
\begin{equation*}
b'_n=\sum_{i=0}^{n-1}(-1)^i d_i.
\end{equation*}
 It turns out that $(C_*(A),b)$ is a chain complex, whose homology is what we call
the \emph{Hochschild homology} of the superalgebra $A$, while the pair $(C_*(A),b')$ is also a
complex, but an acyclic one. To get to the \emph{cyclic} homology of the superalgebra $A$, we need
one more object, namely the \emph{cyclicity map} $t_n:C_n(A)\to C_n(A)$, given by
\begin{equation}\label{cycl}
t_n(a_0\otimes \dots \otimes a_n)=(-1)^{|a_n|(|a_0|+\dots|a_{n-1}|)}a_n\otimes a_0\otimes\dots
\otimes a_{n-1}.
\end{equation}
Finally, we need a last ingredient to define the bicomplex whose total homology is \emph{the cyclic
homology} of the superalgebra $A$, the so-called \emph{norm operator}, $N_n:C_n(A)\to C_n(A)$,
defined by
\begin{equation*}
N_n=1+t_n+\left(t_n\right)^2+\dots+\left(t_n\right)^n.
\end{equation*}
Now, the cyclic homology of the superalgebra $A$ is the homology of the total complex associated to
the following double complex (see \cite{kassel}):
\begin{equation*}
\xymatrix{%
&&&\vdots\ar[d]^{b_3}&\vdots\ar[d]^{b'_3}&\vdots\ar[d]^{b_3}&&\\
j=2&&\cdots&A^{\otimes 3}\ar[l]_{N_2}\ar[d]^{b_2}&A^{\otimes
3}\ar[l]_{1-t_2}\ar[d]^{b'_2}&A^{\otimes
3}\ar[l]_{N_2}\ar[d]^{b_2}&\cdots\ar[l]_{1-t_2}\\
j=1&&\cdots&A^{\otimes 2}\ar[l]_{N_1}\ar[d]^{b_1}&A^{\otimes
2}\ar[l]_{1-t_1}\ar[d]^{b'_1}&A^{\otimes
2}\ar[l]_{N_1}\ar[d]^{b_1}&\cdots\ar[l]_{1-t_1}\\
j=0&&\cdots&A\ar[l]_{N_0}&A\ar[l]_{1-t_0}&A\ar[l]_{N_0}&\cdots\ar[l]_{1-t_0}\\
}
\end{equation*}
Here the even numbered columns are copies of the Hochschild complex, while the odd one are acyclic
complexes.

 There is a close, functorial, relation between the Hochschild and the cyclic homology of  a
superalgebra. This is given by the \emph{Connes' exact sequence} obtained in the following way. Let
us denote by $B_n:C_n(A)\to C_{n+1}(A)$ the operator given by $B_n=(1-t_{n+1})\circ s_n\circ N_n$,
by $S$ the map from the cyclic bicomplex to itself that shifts everything ttwo columns to the left
(\emph{Connes periodicity operator)} and by $I$ the inclusion of the Hochschild complex into the
cyclic bicomplex as the 0-th column. All these maps are compatible with the differentials and
induce maps in homology. Let us use the same notation for these induced maps. Then we have the
following long exact sequence (\cite{kassel}).
\begin{equation}\label{connes}
\begin{split}
\cdots
\xrightarrow{S}HC_{n-1}(A)\xrightarrow{B}HH_n(A)\xrightarrow{I}HC_n(A)\xrightarrow{S}HC_{n-1}(A)
\xrightarrow{B}HH_{n-1}(A)\rightarrow \cdots
\end{split}
\end{equation}
Here $I$ is always surjective and is an isomorphism for $n=0$.
%\begin{thebibliography}{99}
%\end{thebibliography}
\section{The cyclic homology of the superalgebra of supermatrices over a superalgebra $A$}
The intention of this section is to prove the following theorem:
\begin{theo}
Let $A$ be a superalgebra with unit, over a graded-commutative superring $k$. Then, for any pair of
natural numbers $(p,q)$, the supertrace map
\begin{equation*}
\str:\mathcal{M}_{p,q}(A)\to A
\end{equation*}
can be extended to a morphism of the two cyclic bicomplexes:
\begin{equation*}
\Str:CC_{*,*}\left(\mathcal{M}_{p,q}(A)\right)\to CC_{*,*}(A),
\end{equation*}
which, in turn, induces an isomorphism in the cyclic homology. Here we denoted by $CC_{*,*}$ the
chain modules in the cyclic bicomplex.
\end{theo}
\begin{proof}
The proof is an adaptation of the proof for the ungraded case (see, for instance,
\cite{rosenberg}).

Let $M$ be a homogeneous supermatrix of type $p|q$ with entries in the superalgebra $A$. Then, as
we saw earlier, $M$ can be written in the block form
\begin{equation}\label{supermatr1}
M=
\begin{pmatrix}
M_{11}&M_{12}\\
M_{21}&M_{22}
\end{pmatrix}.
\end{equation}
Then the \emph{supertrace} is defined (see \cite{bruzzo},\cite{manin}) as
\begin{equation}\label{supertr}
\str M=\tr M_{11}+(-1)^{1+|M|}\tr M_{22}.
\end{equation}
In a previous paper (\cite{blaga}), we extended this map to a map $\Str$ between the Hochschild
complexes of the two superalgebras and we proved that it induces an isomorphism between the
Hochschild homologies. Namely, we defined $\Str$ on (tensor products of) homogeneous supermatrices
as
\begin{equation}\label{deniss}
\begin{split}
\Str_n\left(M^0\otimes M^1\otimes\dots\otimes M^n\right)=\sum_{i_0=1}^n\sum M^0_{i_0i_1}\otimes
M^1_{i_1i_2}\otimes \dots \otimes M^{n-1}_{i_{n-1}i_n}\otimes M^n_{i_ni_0}+\\
+(-1)^{1+|M^0|+\dots +|M^n|}\sum_{i_0=1}^{p+q}\sum M^0_{i_0i_1}\otimes M^1_{i_1i_2}\otimes \dots
\otimes M^{n-1}_{i_{n-1}i_n}\otimes M^n_{i_ni_0}.
\end{split}
\end{equation}
Here the lower indices of the matrices are not related to the block decomposition, they simple
identify the entries of the matrices. The second sum, in both terms, is taken after all values of
the indices $i_1,\dots, i_n$ (between 1 and $p+q$, of course!).

It turns out that the generalized supertrace is uniquely defined by its values on (tensor products
of) all elementary matrices $E^{ij}(a)$, where $a\in A$ is a homogeneous element  and
\begin{equation}\label{elematr}
E^{ij}_{rs}(a)=
\begin{cases}
a&\text{if}\quad (r,s)=(i,j)\\
0&\text{in all other cases}
\end{cases}.
\end{equation}
Clearly, the parity of $E^{ij}(a)$ is equal to the parity of $a$ if the only non-vanishing element
belongs to a diagonal block and it is equal to the opposite of the parity of $a$ in the rest of the
cases.

We proved (\cite{blaga}) that, on elementary matrices, $\Str$ is easy to compute and we have
\begin{equation}\label{elemsuper}
\Str_n\left(E^{i_0i_1}(a_0)\otimes\dots\otimes E^{i_ni_0}(a_n)\right)=
\begin{cases}
a_0\otimes\dots\otimes a_n, & i_0\leq p,\\
(-1)^{1+\sum\limits_{i=0}^n|a_i|}a_0\otimes\dots\otimes a_n,& i_0>p.
\end{cases}
\end{equation}
Notice that we took a special combination of the indices of the matrices in the argument of $\Str$.
It turns out (and it is very easily checked), that for any other combination of indices the
generalized supertrace vanishes and, as such, it is enough to work with this kind of combinations.

To prove that $\Str$ is a morphism of the cyclic bicomplexes, all we have to prove is to show that
it commutes with the differentials (both horizontal and vertical). But we have already proved that
it commutes with the \emph{vertical} differentials (as it commutes with the faces in the Hochschild
complex (\cite{blaga})). Thus, clearly, all that remains to be seen is that the generalized
supertrace commutes with the cyclicity operator.

Combining the formulas~(\ref{cycl}) and~(\ref{elemsuper}), we get:
\begin{equation}\label{comp1}
\begin{split}
&\left(t_n\circ\Str_n\right)\left(E^{i_0i_1}(a_0)\otimes\dots\otimes E^{i_ni_0}(a_n)\right)=\\
&=
\begin{cases}
(-1)^{|a_n|\sum\limits_{i=0}^{n-1}|a_i|}a_n\otimes a_0\otimes\dots\otimes a_{n-1} &\text{if}\quad
i_0\leq p\\
(-1)^{1+\sum\limits_{i=0}^n|a_i|+|a_n|\sum\limits_{i=0}^{n-1}|a_i|}a_n\otimes
a_0\otimes\dots\otimes a_{n-1} &\text{if}\quad i_0> p
\end{cases}
\end{split}
\end{equation}
and, on the other hand,
\begin{equation}\label{comp2}
\begin{split}
&\left(\Str_n\circ t_n\right)\left(E^{i_0i_1}(a_0)\otimes\dots\otimes E^{i_ni_0}(a_n)\right)=\\
&=
\begin{cases}
(-1)^{\left|E^{i_ni_0}(a_n)\right|\cdot\left|E^{i_0i_n}(a_0\cdots
a_{n-1})\right|}a_n\otimes\dots\otimes a_{n-1}&\text{if}\quad i_n\leq p\\
(-1)^{\left|E^{i_ni_0}(a_n)\right|\cdot\left|E^{i_0i_n}(a_0\cdots
a_{n-1})\right|+1+\sum\limits_{i=0}^{n}|a_i|}a_n\otimes\dots\otimes a_{n-1}&\text{if}\quad i_n>p
\end{cases}.
\end{split}
\end{equation}
Let us denote, for simplicity, $E^{n}(a)=E^{i_0i_1}(a_0)\otimes \dots \otimes E^{i_ni_0}(a_n)$. We
have four situations to consider:
\begin{enumerate}[(i)]
\item $i_0\leq p, i_n\leq p$. In this case, we have
\begin{equation}\label{comp31}
\left(t_n\circ \Str_n\right)(E^n(a))=(-1)^{|a_n|\sum\limits_{i=0}^{n-1}|a_i|}a_n\otimes
a_0\otimes\dots\otimes a_{n-1},
\end{equation}
while
\begin{equation}\label{comp31a}
\begin{split}
\left(\Str_n\circ
t_n\right)(E^n(a))&=(-1)^{\left|E^{i_ni_0}(a_n)\right|\cdot\left|E^{i_0i_n}(a_0\cdots
a_{n-1})\right|}a_n\otimes\dots\otimes a_{n-1}=\\
&=(-1)^{|a_n|\sum\limits_{i=0}^{n-1}|a_i|}a_n\otimes a_0\otimes\dots\otimes a_{n-1}=\left(t_n\circ
\Str_n\right)(E^n(a)).
\end{split}
\end{equation}
\item $i_0\leq p, i_n>p$. In this case, $\left(t_n\circ \Str_n\right)(E^n(a))$ is still given
by~(\ref{comp31}), while
\begin{equation}\label{comp32a}
\begin{split}
\left(\Str_n\circ
t_n\right)(E^n(a))&=(-1)^{(1+|a_n|)\bigg(1+\sum\limits_{i=0}^{n-1}|a_i|\bigg)+1+\sum\limits_{i=0}^{n}|a_i|}a_n\otimes
a_0\otimes\dots\otimes a_{n-1}.
\end{split}
\end{equation}
Having in mind that the parities are elements of $\mathbb{Z}_2$,~(\ref{comp32a}) leads again to
\begin{equation*}
\left(\Str_n\circ t_n\right)(E^n(a))=\left(t_n\circ\Str_n\right)(E^n(a)).
\end{equation*}
\item $i_0>p, i_n\leq p$. We have, now,
\begin{equation}\label{comp33}
\left(t_n\circ
\Str_n\right)(E^n(a))=(-1)^{1+\sum\limits_{i=0}^n|a_i|+|a_n|\sum\limits_{i=0}^{n-1}|a_i|}a_n\otimes
a_0\otimes\dots\otimes a_{n-1}
\end{equation}
and
\begin{equation}\label{comp33a}
\begin{split}
\left(\Str_n\circ t_n\right)(E^n(a))&=
(-1)^{(1+|a_n|)\bigg(1+\sum\limits_{i=0}^{n-1}|a_i|\bigg)}a_n\otimes\dots\otimes
a_{n-1}=\\
&=\left(t_n\circ \Str_n\right)(E^n(a))
\end{split}
\end{equation}
\item $i_0>p, i_n>p$. In this situation, $\left(t_n\circ \Str_n\right)(E^n(a))$ is given
by~(\ref{comp33}) and
\begin{equation}\label{comp34a}
\begin{split}
\left(\Str_n\circ
t_n\right)(E^n(a))&=(-1)^{|a_n|\sum\limits_{i=0}^n|a_i|+1+\sum\limits_{i=0}^n|a_i|}a_n\otimes
a_0\otimes \dots\otimes a_{n-1}=\\
&=\left(t_n\circ \Str_n\right)(E^n(a)).
\end{split}
\end{equation}
\end{enumerate}
We are now, left with the problem of showing that $\Str$ is a quasi-isomorphism between the cyclic
bicomplexes. We are going to do that exactly as in the classical, ungraded case. Namely, by
induction. Since all the cyclic homology groups of negative degree are zero for both superalgebras,
we shall start by $n=0$. But, as we mentioned earlier, as a consequence of the existence of the
Connes long exact sequence, at this level the cyclic and Hochschild homologies coincide, in other
words, we have $HH_0(A)=HC_0(A)$ and the same for the case of the superalgebra
$\mathcal{M}_{p,q}(A)$. Thus, for the zeroth order homology, its is, actually, enough to prove that
$\Str$ is a quasi-isomorphism in the \emph{Hochschild} homology, which we now already for our
previous paper. Let us assume, now, that $\Str$ induces isomorphisms between the cyclic homology
groups up to the degree $n_0-1$, where $n_0$ is a natural number and let us prove that, in this
situation, it also induces an isomorphism in the homology of degree $n_0$. Writing down the Connes
long exact sequences both for $\mathcal{M}_{p,q}(A)$ and $A$, we get the following diagram:
\begin{equation*}
\xymatrix{%
HC_{n_0-1}\big(\mathcal{M}_{p,q}(A)\big)\ar[d]_B\ar[r]&HC_{n_0-1}(A)\ar[d]_B\\
HH_{n_0}\big(\mathcal{M}_{p,q}(A)\big)\ar[d]_I\ar[r]&HH_{n_0}(A)\ar[d]_I\\
HC_{n_0}\big(\mathcal{M}_{p,q}(A)\big)\ar[d]_S\ar[r]&HC_{n_0}(A)\ar[d]_S\\
HC_{n_0-2}\big(\mathcal{M}_{p,q}(A)\big)\ar[d]_B\ar[r]&HC_{n_0-2}(A)\ar[d]_B\\
HH_{n_0-1}\big(\mathcal{M}_{p,q}(A)\big)\ar[r]&HH_{n_0-1}(A)
}
\end{equation*}
The columns of this diagram are exact and the horizontal maps are those induced by the generalized
supertrace in the Hochschild and cyclic homology, respectively. Using the induction hypothesis, as
well as the fact that the generalized supertrace is a quasi-isomorphism at the level of the
Hochschild complexes, it follows that the first two and the last two horizontal maps are
isomorphisms. By the lemma of five (\cite{cartan}) it follows, therefore, that the fifth horizontal
map (the one in the middle) is, also, an isomorphism and, thus, we proved that $\Str$ induces an
isomorphism in the cyclic homology of degree $n_0$, which concludes our proof that $\Str$ is a
quasi-isomorphism at the level of the cyclic bicomplex, as well.
\end{proof}
\section{Final remarks}
The notion of Morita invariance can be extended to the case of graded algebras (see \cite{marcus})
and, in particular, to the case of superalgebras (\cite{blaga1}). It can be shown easily that the
superalgebra of supermatrices of a given type over $A$ and the superalgebra $A$ itself are Morita
equivalent. In the ungraded case, there is well known that the Hochschild and cyclic homology are
Morita invariant. In  \cite{blaga1} we shown (by using a spectral sequences argument) that the same
is true for superalgebras, in the case of Hochschild homology. It is, probably, possible to adapt
the arguments of McCarthy (\cite{mccarthy}) to the ``supercase'', obtaining, in this way, a proof
of the Morita invariance of the cyclic homology for superalgebras.

The advantage of the approach we used here is that (in the particular case of superalgebras of
supermatrices), we can provide explicitly the isomorphism, using some elementary constructions.

\end{document}